\documentclass[12pt]{article}

\usepackage{amsthm,amsmath,amssymb,fullpage,mathrsfs}

\title{{On the characterization of graphs with tree 3-spanners}}
\author{{Lan Lin$^{a}$\thanks{Corresponding Author. E-mail address:
linlan@tongji.edu.cn}, \ \ \ Yixun Lin$^{b}$}\\
{\small $^{a}$ School of Electronics and Information Engineering,
Tongji University,}\\ {\small Shanghai 200092, China}\\
{\small $^{b}$ School of Mathematics and Statistics, Zhengzhou University,}\\
{\small Zhengzhou 450001, China}}

\date{}
\begin{document}
\maketitle
\renewcommand\baselinestretch{1.0}

{\indent\bf Abstract.} \ The tree spanner problem for a graph $G$ is
as follows: For a given integer $k$, is there a spanning tree $T$ of
$G$ (called a tree $k$-spanner) such that the distance in $T$ between
every pair of vertices is at most $k$ times their distance in $G$?
The minimum $k$ that $G$ admits a tree $k$-spanner is denoted by $\sigma(G)$.
It is well known in the literature that determining $\sigma(G)\leq 2$
is polynomially solvable, while determining $\sigma(G)\leq k$ for
$k\geq 4$ is NP-complete. A long-standing open problem is to characterize
graphs with $\sigma(G)=3$. This paper settles this open problem
by proving that it is polynomially solvable.

{\indent\bf Keywords.} \ spanning tree optimization, tree spanner,
characterization, polynomial-time algorithm

\section{Introduction}
\hspace*{0.5cm} The tree spanner problem originally arises from Peleg and
Ullman \cite{Peleg89} in a decision version: For a given integer $k$, is
there a spanning tree $T$ of $G$ (called a tree $k$-spanner) such that the
distance in $T$ between every pair of vertices is at most $k$ times their
distance in $G$? We prefer to state this problem in an optimization version.

Let $G$ be a simple connected graph with vertex set $V(G)$ and edge set $E(G)$.
A spanning tree $T$ of $G$ is a subgraph of $G$ which is a tree with $V(T)=V(G)$.
For each edge $uv\in E(G)$, denote by $d_T(u,v)$  the distance between $u$ and $v$
in $T$, that is, the length of the unique $u$-$v$-path in $T$. The maximum distance
in $T$ between two adjacent vertices
\begin{equation}
\sigma(G,T):=\max_{uv\in E(G)}\,d_T(u,v)
\end{equation}
is called {\it the max-stretch} of spanning tree $T$. {\it The minimum stretch
spanning tree problem} (or the MSST problem) is to find a spanning tree $T$ such
that $\sigma(G,T)$ is minimized, where the minimum value
\begin{equation}
\sigma(G):=\min \{\sigma(G,T): T \mbox{ is a spanning tree of } G\}
\end{equation}
is called the {\it tree-stretch} of $G$ (see survey \cite{Lieb08}). A
spanning tree $T$ attaining this minimum value is called an {\it optimal
spanning tree}. To compare with the above decision version, a spanning
tree $T$ is a tree $k$-spanner if and only if $\sigma(G,T)\leq k$, that
is, $d_T(u,v)\leq k$ for every edge $uv\in E(G)$ (see Observations in
\cite{Cai95}).

In the algorithmic aspect, the MSST problem has been proved NP-hard
\cite{Brand04,Cai95,Fekete01}, and the fixed-parameter polynomial algorithms
have been presented in \cite{Brand04,Fekete01}. For the characterization problems,
it is known that determining $\sigma(G)\leq 2$ is polynomially solvable
(\cite{Bondy89,Cai95}), while determining $\sigma(G)\leq k$ for $k\geq 4$ is
NP-complete (\cite{Cai95} among others). The problem of characterizing $\sigma(G)=3$
was a remarkable open problem in this area (it was conjectured NP-complete
in \cite{Cai95,Fekete01}).

Some partial results have been obtained in the literature. For example,
Madanlal et al. \cite{Madan96} showed that $\sigma(G)\leq 3$ for
all interval and permutation graphs, and that a regular bipartite
graph $G$ has $\sigma(G)\leq 3$ if and only if it is complete.
Brandst\"adt et al. \cite{Brand04} and Venkatesan et al. \cite{Venka97}
showed that $\sigma(G)\leq 3$ for interval graphs and split graphs $G$.
Moreover, Fekete and Kremer \cite{Fekete01} showed that $\sigma(G)\leq 3$
can be decided in polynomial time for planar graphs. The characterization
of directed path graphs with $\sigma(G)\leq 3$ was studied in \cite{Le99,
Panda07,Panda09}. Some exact results for special graphs can be seen in
\cite{Lin20,Lin21}.

In this paper, we establish a necessary and sufficient condition for
graphs with $\sigma(G)\leq 3$ and present the corresponding
recognition algorithm, which has polynomial-time complexity.
Consequently, we conclude that the problem of characterizing
$\sigma(G)\leq 3$ is polynomially solvable.

The paper is organized as follows. In Section 2, some definitions
and notations are introduced. In Section 3, we prove the necessary and
sufficient condition of optimal spanning trees. Section 4 is devoted to
the polynomial-time recognition algorithm for $\sigma(G)\leq 3$.
Section 5 contains a brief summary.

\section{Preliminaries}
\hspace*{0.3cm} We follow the graph-theoretic terminology and notation of
\cite{Bondy08}. For a subset $S\subseteq V(G)$, the {\it neighbor set} of
$S$ is defined by $N_G(S):=\{v\in V(G)\setminus S:\exists u\in S,uv\in E(G)\}$.
We abbreviate $N_G(\{v\})$ to $N_G(v)$ for a vertex $v\in V(G)$. Furthermore,
we denote the {\it closed neighbor set} of $v\in V(G)$ by $N_G[v]:=\{v\}
\cup N_G(v)$. For $S\subseteq V(G)$, we denote by $G[S]$ the subgraph
induced by $S$, which is formed by $S$ and all edges with both ends in
$S$.

Let $T$ be a spanning tree of $G$. Usually, the spanning tree $T$ is
also regarded as a set of edges. The {\it cotree} $\overline T$ of $T$
is the complement of $T$ in $E(G)$, namely $\overline T=E(G)\setminus T$.
For each $e\in \overline T$, the unique cycle in $T+e$ is a {\it fundamental
cycle} with respect to the cotree edge $e$. So, the MSST problem is
equivalent to a problem of finding a spanning tree $T$ such that the
length of the longest fundamental cycle is minimized, where the tree-stretch
$\sigma(G)$ is one less than this minimum length of longest fundamental cycle.

In a duality point of view, for each $e\in T$, the edge cut between two components
of $T-e$ is a {\it fundamental edge cut} (cocycle). In more detail, let
$X_e$ and $Y_e$ be the vertex sets of two components of $T-e$. Write
$\partial (X_e):=\{uv\in E(G):u\in X_e, v\in Y_e\}$. Then $\partial (X_e)$
is the fundamental edge cut with respect to $e$. A prominent property is that
for $e'\in \overline T$, $e$ is contained in the fundamental cycle with
respect to $e'$ if and only if $e'$ is contained in the fundamental
edge cut $\partial (X_e)$ with respect to $e$ (cf. \cite{Bondy08}).

In a spanning tree $T$ of $G$, a vertex $v$ is called a {\it pendant vertex}
of $T$ if it is incident with only one edge in $T$ (it is also called a {\it leaf}
of $T$). An edge $uv\in T$ is called an {\it outer edge} of $T$ if either $u$
or $v$ is a pendant vertex in $T$. An edge $uv\in T$ is called an {\it inner
edge} of $T$ if neither $u$ nor $v$ is pendant.

A {\it star} is a tree with diameter two, in which the only non-pendant
vertex is called its {\it center}. In a star $T$, all edges are outer.
A {\it double star} is a tree with diameter three, that is, it can be
obtained from two stars by identifying one of their edges, where the
common edge of two stars is called the {\it central edge}. In a double
star $T$, all edges are outer except the central edge (which is inner).

The MSST problem is closely related to the minimum diameter spanning tree
problem, which is to find a spanning tree with minimum diameter. Here,
the minimum diameter is denoted by
\begin{equation}
D_T(G):=\min \{\max_{u,v\in V(G)}\,d_T(u,v): T \mbox{ is a spanning
tree of } G\}.
\end{equation}
It is known in Hassin and Tamir \cite{Hassin95} that the minimum diameter
spanning tree problem can be easily solved by the absolute 1-center problem
in $O(mn+n^2\log n)$ time (where $n$ is the number of vertices and $m$ the
number of edges of $G$) as follows. Let $A(G)$ denote the continuum set of
points on the edges of $G$. In the absolute 1-center problem, the point
$x\in A(G)$ is continuously located on the edges of $G$ ($x$ can be
regarded as a new vertex subdividing some edge). A point $x^*\in A(G)$
is an absolute 1-center of $G$ if the function
$$ F(x):=\max_{v\in V(G)}d_G(x,v) \quad (x\in A(G))$$
attains its minimum at $x^*$. For this problem, there has been simple
algorithms in graph theory and location problems (see, e.g., \cite{Kariv79}).
Let $x^*$ be an absolute 1-center of $G$. Then a shortest path tree $T(x^*)$
connecting $x^*$ to all vertices in $V(G)$ is a minimum diameter spanning
tree of $G$. This shortest path tree can be constructed in $O(n^2)$ time
by running the shortest path algorithm. By comparing (1)(2) and (3), we
obtain the following upper bound of the tree-stretch:
$$\sigma(G)\leq D_T(G).$$

A {\it block} of $G$ is a subgraph of $G$ which contains no cut
vertices and it is maximal with respect to this property. Two blocks
of $G$ have at most one vertex (a cut vertex) in common. As each
fundamental cycle is contained in a block, we have the following
observation: If $G$ has blocks $G_1,G_2,\ldots,G_k$, then
$\sigma(G)=\max_{1\leq i\leq k}\sigma(G_i)$. Hence $\sigma(G)\leq 3$
if and only if $\sigma(G_i)\leq 3$ for every block $G_i$ ($1\leq i\leq k$).
So, we may assume that $G$ is itself a block, that is a 2-connected graph
(for $n\geq 3$). Consequently, we need only discuss the characterization
problem for 2-connected graphs henceforth.

The {\it girth} of $G$ is the length of a shortest cycle in $G$, denoted
by $g(G)$. Since the length of a fundamental cycle is at least the girth
$g(G)$, we have a lower bound that
$$\sigma(G)\geq g(G)-1.$$
It is interesting to characterize the graphs with the minimum value of
$\sigma(G)=g(G)-1$. For $g(G)=3$, this characterization problem has been
settled, as stated below.

Clearly, $\sigma(G)=1$ if and only if $G$ is a tree. A graph $G$ with
$\sigma(G)=2$ is called a {\it trigraph} in \cite{Bondy89}, that is,
it has a spanning tree $T$ for which every fundamental cycle is a
triangle. For this class of graphs, Bondy \cite{Bondy89} and Cai and
Corneil \cite{Cai95} presented the following characterization. Herein,
a vertex of $G$ is called a {\it dominating vertex} (or a {\it universal
vertex}) if it is adjacent to all other vertices of $G$. A {\it 3-connected
component} of a 2-connected graph $G$ is a maximal subgraph of $G$ that
has no 2-vertex cut. Note that a 3-connected component with at least four
vertices is 3-connected (apart from the trivial case $K_3$).

{\bf Theorem 2.1} (\cite{Bondy89,Cai95}).  Suppose that $G$ is a 2-connected
graph. Then $G$ admits $\sigma(G)=2$ with an optimal spanning tree $T$ if and
only if \\
{\indent} (i) if $G$ has no 2-vertex cuts, then $G$ contains a dominating
vertex, and $T$ is a spanning star of $G$; or \\
{\indent} (ii) if $G$ has a 2-vertex cut $\{u,v\}$, then $uv\in T$ and for
each 3-connected component $H$ containing $\{u,v\}$ in $G$, the subtree
$T\cap H$ is a spanning star of $H$ with center $u$ or $v$.

In the sequel, we aim to generalize Theorem 2.1 to be a characterization
of $\sigma(G)\leq 3$.

\section{Characterization of optimal spanning trees}
\hspace*{0.5cm} Let us take a closer look on the characterization of
$\sigma(G)=2$. The following necessary and sufficient condition is equivalent
to Theorem 2.1. Let $S$ be a vertex cut of $G$. Suppose that $V_1,V_2,\ldots,V_k$
are the vertex sets of the components of $G-S$. Then the subgraph $G[S\cup V_i]$
induced by $S\cup V_i$ is called an {\it $S$-component} of $G$ (see Section 9.4 of
\cite{Bondy08}).

{\bf Theorem 3.1.}  Suppose that $G$ is a 2-connected graph. Then $G$
admits $\sigma(G)=2$ with an optimal spanning tree $T$ if and only if \\
{\indent} (a) $T$ has diameter two ($T$ is a star with all edges
being outer); or \\
{\indent} (b) for each inner edge $uv$ of $T$, the set $\{u,v\}$ is a 2-vertex
cut of $G$ and for each $\{u,v\}$-component $H$, the edge $uv$ is outer in $T\cap H$.

{\bf Proof.}  To show sufficiency, suppose that there exists a spanning
tree $T$ satisfying (a) or (b). If (a) holds, namely, $T$ is a spanning
star, then for every cotree edge $e\in \overline T$, the fundamental cycle
with respect to $e$ is a triangle, and thus $\sigma(G)=2$. Next, we show that
if (b) holds, then $d_T(x,y)\leq 2$ for every cotree edge $xy\in \overline T$.
Assume, to the contrary, that $d_T(x,y)\geq 3$ for some $xy\in \overline T$. Then
the path connecting $x$ and $y$ in $T$ is in the form of $x\cdots x'uvy'\cdots y$,
where $x'\in N_T(u)$ and $y'\in N_T(v)$ (maybe $x'=x$ and/or $y'=y$). Hence $uv$
is an inner edge of $T$ and it follows from (b) that $\{u,v\}$ is a 2-vertex cut
of $G$. However, in view of $xy\in E(G-\{u,v\})$, we see that the path $ux'\cdots
xy\cdots y'v$ is contained in a $\{u,v\}$-component $H$ of $G$. As a result, $uv$
is not outer in $T\cap H$, contradicting (b). Therefore, we show that $d_T(x,y)\leq 2$
for every cotree edge $xy\in \overline T$, thus $\sigma(G)=2$.

To show necessity, suppose that $\sigma(G)=2$ and $T$ is an optimal spanning
tree. By Theorem 2.1, if $G$ has no 2-vertex cuts, then (i) implies (a).
Otherwise $G$ has 2-vertex cuts. We proceed to show assertion (b). In fact,
if $uv$ is an inner edge of $T$ but $\{u,v\}$ is not a 2-vertex cut of $G$,
then $G-\{u,v\}$ is connected. For $uv\in T$, we can consider the
fundamental edge cut $\partial (X_{uv})$ between two components of $T-uv$,
where $X_{uv}$ and $Y_{uv}$ are the vertex sets of the components of $T-uv$
containing $u$ and $v$ respectively. As $G-\{u,v\}$ is connected, there
must be an edge $xy\in E(G)$ such that $\{x,y\}\cap \{u,v\}=\emptyset$ and
$xy\in \partial (X_{uv})$. Thus $xy\in \overline T$. Since the path connecting
$x$ and $y$ in $T$ contains the inner edge $uv$, it follows that $d_T(x,y)\geq 3$,
which contradicts $\sigma(G)=2$. Thus we show that $\{u,v\}$ is indeed a 2-vertex
cut of $G$. Furthermore, we can show that for each $\{u,v\}$-component $H$, the
edge $uv$ is outer in $T\cap H$. Suppose the contrary. Then there exists a
$\{u,v\}$-component $H$ such that $uv$ is not outer in $T\cap H$ (i.e., neither
$u$ nor $v$ is pendant in  $T\cap H$), then there must be $x\in N_T(u)\setminus
\{v\}$ and $y\in N_T(v)\setminus \{u\}$ contained in $T\cap H$. Hence $x\in X_{uv}$
and $y\in Y_{uv}$. In this $\{u,v\}$-component $H=G[\{u,v\}\cup V_i]$, the path
connecting $x$ and $y$ in $V_i$ must meet the fundamental edge cut $\partial (X_{uv})$
at some edge $x'y'\in \partial (X_{uv})$ (maybe $x'=x$ and/or $y'=y$). Thus
$x'y'\in \overline T$. Since the path connecting $x'$ and $y'$ in $T$ contains
$uv$, it follows that $d_T(x',y')\geq 3$, which contradicts $\sigma(G)=2$.
Therefore, (b) is proved. The proof is complete. $\Box$

Now we start considering the characterization of $\sigma(G)\leq 3$.

For a spanning tree $T$ and a non-pendant vertex $v\in V(T)$, the subset
$S=\{v\}\cup N_T(v)$ is called a {\it star-set} of $T$ with center $v$,
which induces a maximal subtree of diameter two in $T$. In this star-set,
a vertex $u\in S\setminus \{v\}$ is called a {\it contact vertex} of $S$
if it is not a pendant vertex of $T$. Moreover, a star-set $S$ is said to
be {\it outer} if it has at most one contact vertex. In contrast to that,
a star-set $S$ is said to be {\it inner} if it has at least two contact
vertices, which are adjacent to some vertices outside $S$. The following
is a generalization of Theorem 3.1.

{\bf Theorem 3.2.}  Suppose that $G$ is a 2-connected graph. Then $G$ admits
$\sigma(G)\leq 3$ with an optimal spanning tree $T$ if and only if \\
{\indent} (I) $T$ has diameter at most three ($T$ is a star or double star,
which has no inner star-sets); or \\
{\indent} (II) each inner star-set $S$ of $T$ is a vertex cut of $G$ and
for each $S$-component $H$ of $G$, the star-set $S$ is outer in $T\cap H$.

{\bf Proof.}  We first prove sufficiency. Suppose that there exists a spanning
tree $T$ satisfying condition (I) or (II). If (I) holds and $T$ is a spanning
star, then $\sigma(G)=2$ as in Theorem 3.1. If $T$ is a double star, then
for every cotree edge $e\in \overline T$, the fundamental cycle with respect
to $e$ is a triangle or quadrilateral, and thus $\sigma(G)\leq 3$. Otherwise,
the diameter of $T$ is greater than three and thus it has inner star-sets (if
the diameter path of $T$ contains $xuvwy$, then the star-set with center $v$
is inner). We next show that if condition (II) holds, then $d_T(x,y)\leq 3$ for
every cotree edge $xy\in \overline T$. Assume, to the contrary, that $d_T(x,y)
\geq 4$ for some $xy\in \overline T$. Without loss of generality we may assume
that $d_T(x,y)=4$ (as the proof is almost the same for longer distance). Then
the path connecting $x$ and $y$ in $T$ is in the form of $xuvwy$, where $x\in
N_T(u)$ and $y\in N_T(w)$. We take a star-set $S=\{v\}\cup N_T(v)$. Since
$u,w\in S$ are not pendant vertices of $T$, they must be contact vertices of $S$.
Hence $S$ is an inner star-set of $T$. By condition (II), $S$ is a vertex cut
of $G$. Since $xy\in \overline T\subseteq E(G)$, we see that the cycle $xuvwyx$
is contained in an $S$-component $H$ containing $N_T(u)$ and $N_T(w)$. Consequently,
the two vertices $u,w\in S$ are also contact vertices in $T\cap H$ (as shown in
Figure 1). Thus $S$ is not outer in $T\cap H$, which contradicts (II). Therefore,
we show that $d_T(x,y)\leq 3$ for every cotree edge $xy\in \overline T$, and so
$\sigma(G)\leq 3$.

\begin{center}
\setlength{\unitlength}{0.5cm}
\begin{picture}(11,11)

\multiput(4,4)(2,0){2}{\circle*{0.3}}
\multiput(3,5)(2,0){2}{\circle*{0.3}}
\multiput(1,8)(5,0){2}{\circle*{0.3}}
\multiput(10,6)(0,0){1}{\circle*{0.3}}

\put(4,4){\line(1,0){2}} \put(6,4){\line(2,1){4}}
\put(4,4){\line(-1,1){1}} \put(4,4){\line(1,1){1}}
\put(3,5){\line(-2,3){2}} \put(5,5){\line(1,3){1}}
\bezier{30}(1,8)(3.5,8)(6,8)
\put(4.5,4.5){\oval(5,3)}
\put(1,8.5){\oval(3,4)}
\put(6,8.5){\oval(3,4)}
\put(10,6.5){\oval(3,4)}

\put(3.3,6.8){\makebox(1,0.5)[l]{\small $H$}}
\put(2.5,3.7){\makebox(1,0.5)[l]{\small $S$}}
\put(1.3,8.1){\makebox(1,0.5)[l]{\small $x$}}
\put(6.2,8.1){\makebox(1,0.5)[l]{\small $y$}}
\put(3.3,5){\makebox(1,0.5)[l]{\small $u$}}
\put(5.3,5){\makebox(1,0.5)[l]{\small $w$}}
\put(3.9,3.2){\makebox(1,0.5)[l]{\small $v$}}
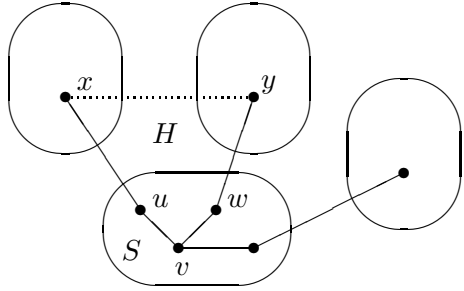
\put(-3,0){\makebox(1,0.5)[l]{\small Figure 1. Illustration for the proof of sufficiency.}}
\end{picture}
\end{center}

We next show necessity. Suppose that $\sigma(G)\leq 3$ and $T$ is an optimal
spanning tree. If $T$ has diameter at most three, then (I) holds. Otherwise,
the diameter of $T$ is at least four. Then $T$ contains inner star-sets. Now
let $S$ be an inner star-set with center $v$ and with at least two contact
vertices $u$ and $w$. We proceed to show assertion (II) by contradiction.
Assume, to the contrary, that $S$ is not a vertex cut of $G$ or it is a vertex
cut but it is not outer in an $S$-component $H$. If $S$ is not a vertex cut of $G$,
then $G-S$ is connected. We can take $x\in N_T(u)\setminus S$ and $y\in N_T(w)
\setminus S$ such that $x$ and $y$ are connected by a path $P$ in $G-S$.
Similarly, if $S$ is a vertex cut of $G$ but it is not outer in an $S$-component
$H$, then $S$ has two contact vertices $u$ and $w$ with respect to $H$. We
can also take $x\in N_T(u)\setminus S$ and $y\in N_T(w) \setminus S$ such that
$x$ and $y$ are connected by a path $P$ in $H$. Note that this path $P$ cannot
be contained in $T$ (for otherwise $P$ and $xuvwy$ form a cycle in $T$). So,
$P$ contains some edges in $\overline T$. We distinguish two cases as follows.\\
{\indent} $\bullet$ The path $P$ contains only one edge in $\overline T$.
Suppose that $P=x\cdots x'y'\cdots y$, where $x'y'\in \overline T$ and
the other edges of $P$ are contained in $T$. Then the path connecting $x'$
and $y'$ in $T$ is $x'\cdots xuvwy\cdots y'$, and thus $d_T(x',y')\geq 4$,
which contradicts that $\sigma(G)\leq 3$. \\
{\indent} $\bullet$ The path $P$ contains more than one edge in
$\overline T$. Clearly, $S$ is a vertex-cut of $T$ that separates
$N_T(u)\setminus S$ and $N_T(w)\setminus S$ in different components
of $T-S$. Let $T_x$ and $T_y$ be the components of $T-S$ containing
$x$ and $y$ respectively. By assumption, the path $P$ contains two or
more edges in $\overline T$. Let $x'z'$ be the first cotree edge in $P$,
where $x'$ belongs to the component $T_x$ of $T-S$ and $z'$ belongs to
another component $T_z$ of $T-S$. Let $z$ be the vertex in $T_z$ which is
adjacent to some contact vertex $t\in N_T(v)$ (where $t$ is not a pendant
vertex of $T$). In this situation, we can consider the contact vertices
$u$ and $t$ (instead of $u$ and $w$) in $S$. Then the path $P'$ connecting
$x\in N_T(u)\setminus S$ and $z\in N_T(t)\setminus S$ contains only one
edge in $\overline T$ (as shown in Figure 2). Then the path connecting $x'$
and $z'$ in $T$ is $x'\cdots xuvtz\cdots z'$, and thus $d_T(x',z')\geq 4$,
which is also a contradiction to $\sigma(G)\leq 3$. To summarize, we show
assertion (II). The proof is complete. $\Box$

\begin{center}
\setlength{\unitlength}{0.5cm}
\begin{picture}(11,11)

\multiput(5,4)(0,0){1}{\circle*{0.3}}
\multiput(3,5)(2,0){3}{\circle*{0.3}}
\multiput(1,7)(4,0){3}{\circle*{0.3}}
\multiput(1,8)(4,0){3}{\circle*{0.3}}

\put(5,4){\line(-2,1){2}} \put(5,4){\line(0,1){1}}
\put(5,4){\line(2,1){2}} \put(3,5){\line(-1,1){2}}
\put(5,5){\line(0,1){3}} \put(7,5){\line(1,1){2}}
\put(1,7){\line(0,1){1}} \put(9,7){\line(0,1){1}}
\put(5,8){\line(1,0){1.7}} \put(8,8){\line(1,0){1}}
\bezier{50}(1,8)(5,8)(9,8)
\put(5,4.1){\oval(6,3)}
\put(1,8.2){\oval(2.7,4)}
\put(5,8.2){\oval(2.7,4)}
\put(9,8.2){\oval(2.7,4)}

\put(0.8,9){\makebox(1,0.5)[l]{\small $T_x$}}
\put(4.8,9){\makebox(1,0.5)[l]{\small $T_z$}}
\put(8.8,9){\makebox(1,0.5)[l]{\small $T_y$}}
\put(3.5,3){\makebox(1,0.5)[l]{\small $S$}}
\put(2.6,8.1){\makebox(1,0.5)[l]{\small $P'$}}
\put(1.3,8.2){\makebox(1,0.5)[l]{\small $x'$}}
\put(1.3,7){\makebox(1,0.5)[l]{\small $x$}}
\put(5.3,8.2){\makebox(1,0.5)[l]{\small $z'$}}
\put(5.3,7){\makebox(1,0.5)[l]{\small $z$}}
\put(9.3,8.2){\makebox(1,0.5)[l]{\small $y'$}}
\put(9.3,7){\makebox(1,0.5)[l]{\small $y$}}

\put(3.3,5){\makebox(1,0.5)[l]{\small $u$}}
\put(5.3,4.8){\makebox(1,0.5)[l]{\small $t$}}
\put(7.1,4.2){\makebox(1,0.5)[l]{\small $w$}}
\put(4.8,3.2){\makebox(1,0.5)[l]{\small $v$}}
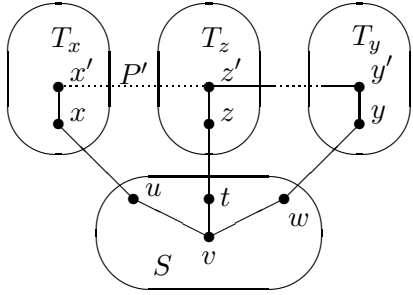
\put(-3,0){\makebox(1,0.5)[l]{\small Figure 2. Illustration for the proof of necessity.}}
\end{picture}
\end{center}

Note that the condition (II) in Theorem 3.2 can be stated in different ways.
For an inner star-set $S=\{v\}\cup N_T(v)$ in $T$, the subtrees of $T-S$
are called {\it $S$-branches} of $T$. For example, $T_x,T_y,T_z$ in Figure 2
are $S$-branches.

{\bf Proposition 3.3.}  Suppose that the inner star-set $S$ of $T$ is a vertex cut of $G$.
The following assertions are equivalent:\\
{\indent} (IIa) For each $S$-component $H$ of $G$, the star-set $S$ is outer in $T\cap H$.\\
{\indent} (IIb) For any two contact vertices $u$ and $w$ of $S$, $N_T(u)\setminus S$ and
$N_T(w)\setminus S$ are contained in different components of $G-S$.\\
{\indent} (IIc) The vertices in different $S$-branches are independent (nonadjacent).

{\bf Proof.} (IIa) $\Rightarrow$ (IIb). If for two contact vertices $u$ and $w$
of $S$, $N_T(u)\setminus S$ and $N_T(w)\setminus S$ are contained in the same component
of $G-S$, then $u$ and $w$ are contained in the same $S$-component $H$, thus $S$ contains
two contact vertices $u$ and $w$ in $H$, and so $S$ is not outer in $T\cap H$.\\
{\indent} (IIb) $\Rightarrow$ (IIc). If two vertices $x$ and $y$ contained in
different $S$-branches leading from the contact vertices $u$ and $w$ respectively
are adjacent, then $N_T(u)\setminus S$ and $N_T(w)\setminus S$ are contained
in the same component of $G-S$.\\
{\indent} (IIc) $\Rightarrow$ (IIa). If $S$ is not outer in an $S$-component $H$, then
$H$ contains two contact vertices $u$ and $w$ of $S$. Thus two different $S$-branches
leading from $u$ and $w$ are contained in the same $S$-component $H$. Hence there must be
two vertices $x$ and $y$ contained in different $S$-branches such that they are adjacent.
This completes the proof. $\Box$

To summarize roughly, for $\sigma(G)=1$, $G$ (a tree) is decomposable by an inner
vertex; For $\sigma(G)=2$, $G$ is decomposable by a 2-vertex set of inner edge;
For $\sigma(G)=3$, $G$ is decomposable by an inner star-set. For the moment, if
we can establish a `non-trivial' necessary and sufficient condition (or a `good'
characterization) for the solutions of an optimization problem, then this problem
is unlikely to be NP-hard.

\section{Recognition algorithm}
\hspace*{0.5cm} By virtue of Theorem 3.2, in order to decide whether $G$
admits $\sigma(G)\leq 3$, it suffices to find a spanning tree $T$ satisfying
condition (I) or (II). To check condition (I), we can apply the minimum
diameter spanning tree algorithm (Hassin and Tamir \cite{Hassin95}), in
which the minimum spanning-tree diameter $D_T(G)$ is determined. As stated
earlier, the running time of this algorithm is $O(mn+n^2\log n)$. However, we
need only use the algorithm to decide whether $D_T(G)\leq 3$ (namely, $T$ is
either a star or a double star). This is to decide whether there exists an edge
$xy$ such that $\{x,y\}$ is a dominating set of $G$. So, we have a simplified
algorithm as follows.

{\noindent}{\bf Basic Procedure} \\
{\noindent}(1) For $xy\in E(G)$, set $T:=\{xy\}$. \\
{\noindent}(2) Let $\{v_1,v_2,\ldots,v_k\}$ denote the set of vertices not covered by $T$.\\
{\bf For} $i:=1$ {\bf to} $k$ {\bf do}:\\
{\indent} {\bf If} $v_ix\in E(G)$, {\bf then} set $T:=T\cup \{v_ix\}$ {\bf else}\\
{\indent} $\quad$ {\bf If} $v_iy\in E(G)$, {\bf then} set $T:=T\cup \{v_iy\}$ {\bf else} go to (4).\\
{\noindent}(3) Return the double star $T$ with central edge $xy$ or the star $T$ with center $x$ or $y$.\\
{\noindent}(4) Stop with answer `no' ($\{x,y\}$ is not a dominating set).

{\bf Lemma 4.1.} \ $D_T(G)\leq 3$ can be determined in $O(mn)$ time.

{\bf Proof.} \ We can run the Basic Procedure for all edges $xy\in E(G)$. If the
procedure yields a spanning tree $T$ that is a star or double star, then $D_T(G)\leq 3$.
Otherwise there exists no such dominating set and thus $D_T(G)>3$. For each edge $xy$,
the procedure is performed in $O(n)$ time and the number of procedures is $m$. Hence
the overall running time is $O(mn)$. $\Box$

In what follows a decomposition approach can be established by using condition (II)
of Theorem 3.2.

{\bf Lemma 4.2.} \ Suppose that $D_T(G)>3$. If $G$ admits $\sigma(G)\leq 3$, then
there exists a vertex $v$ such that $N_G[v]=\{v\}\cup N_G(v)$ is a vertex cut of $G$.

{\bf Proof.} \ Suppose that $D_T(G)>3$ and $\sigma(G)\leq 3$. Let $T$ be an optimal
spanning tree with minimum diameter $D(T)$. Then $D(T)\geq D_T(G)>3$ and so $T$
contains inner star-sets. Let $S=\{v\}\cup N_T(v)$ be an inner star-set. By (II)
of Theorem 3.2, $S$ is a vertex cut of $G$. We proceed to show that $N_G[v]=\{v\}
\cup N_G(v)$ is a vertex cut of $G$. As $N_T(v)\subseteq N_G(v)$, we can move the
vertices of $N_G(v)\setminus N_T(v)$ into $S$ so that $N_G[v]$ is a bigger vertex cut.
This can be done (thus the proof is complete) if $G-N_G[v]$ has at least two components.
Otherwise $G-N_G[v]$ has only one component, i.e., $G-N_G[v]$ is connected. Let
$V_1,V_2,\ldots,V_k$ be the vertex sets of the components of $G-S$ ($k\geq 2$).
Then all $V_1,V_2,\ldots,V_k$ but one, say $V_k$, are contained in $N_G[v]$. That
is to say, all vertices in $V_1\cup \cdots \cup V_{k-1}$ are adjacent to $v$ in $G$.
We can construct another spanning tree $T'$ from $T$ by replacing each path in $T$
from $v$ to $z\in V_1\cup \cdots\cup V_{k-1}$ by a pendant edge $vz$. For this spanning
tree $T'$, any cotree-edge incident with vertices of $V_1\cup \cdots\cup V_{k-1}$ has
stretch two and so $\sigma(G,T')\leq 3$ as before. Thus $T'$ is also an optimal spanning
tree of $G$. However, the diameter $D(T')$ is $\max\{d_{T'}(v,x):x\in V_k\}+1$, which is
less than $D(T)$ (since $S$ is an inner star-set, $\max\{d_{T}(v,x):x\in V_i\}>1$ for at
least one $i$ with $1\leq i\leq k-1$). This yields a contradiction to the minimality of
$D(T)$. The proof is complete. $\Box$

If there exists no vertex $v$ such that $N_G[v]$ is a vertex cut of $G$, then
$D_T(G)>3$ implies $\sigma(G)>3$ by Lemma 4.2 and the decision process is finished.
In other words, if $G$ has no vertex cuts $N_G[v]$ and $G$ admits $\sigma(G)\leq 3$,
then it holds that $D_T(G)\leq 3$, which can be decided in polynomial time by Lemma 4.1.

This kind of vertex cut $X=N_G[v]$ may be called a {\it neighbor cut} of $G$. In the proof
of Lemma 4.2, it is easy to get a neighbor cut $X=N_G[v]$ from an inner star-set $S=
N_T[v]$ by simply adding some vertices. The reverse problem, that is, to find out the inner
star-set contained in a given neighbor cut, is more difficult. This is what we try hard
to solve later. Note that we can determine all neighbor cuts $X=N_G[v]$ in $O(n^2)$ time.
This is because for each vertex $v\in V(G)$, we can check whether $N_G[v]$ is a vertex cut
by computing the number of components of $G-N_G[v]$ in $O(n)$ time and this procedure is
applied at most $n$ time. So, we can assume that we have a neighbor cut in hand.

Suppose that $X$ is a neighbor cut of $G$, and $V_1,V_2,\ldots,V_k$ are the vertex sets of
the components of $G-X$ ($k\geq 2$). As before, the subgraph $G_i:=G[X\cup V_i]$ induced by
$X\cup V_i$ is called an {\it $X$-component} of $G$ ($1\leq i\leq k$). We say that $G$ is
decomposed by $X$ into $G_1,G_2,\ldots,G_k$. We have a basic decomposability approach as follows.

{\bf Lemma 4.3.} \ Suppose that $X=N_G[v]$ is a neighbor cut of $G$ and $G_i:=G[X\cup V_i]$
($1\leq i\leq k$) are decomposed by $X$. Then $\sigma(G)\leq 3$ if and only if $\sigma(G_i)\leq 3$
for all $X$-components $G_i$ ($1\leq i\leq k$) and the optimal spanning trees $T_i$ of $G_i$
are coincident in $G[X]$.

{\bf Proof.} \ For any spanning tree $T$ of $G$, the subtree $T_i$ restricted in $G_i$
is a spanning tree of $G_i$. For any cotree edge $uv\in \overline T$ in $G$, it must be
contained in some $G_i$; if it is contained in $G_i$, then it is a cotree edge with respect
to $T_i$. If $\sigma(G)\leq 3$, then for an optimal spanning tree $T$, the subtree $T_i$
restricted in $G_i$ satisfies $\sigma(G_i,T_i)\leq 3$, and all subtrees $T_i$ are
coincident with $T\cap G[X]$.

Conversely, if for all $G_i$ with $1\leq i\leq k$, we have $\sigma(G_i)\leq 3$ and their
optimal spanning trees $T_i$ are coincident in $G[X]$, then we can construct a spanning tree
$T$ of $G$ by connecting $T_1,T_2,\ldots,T_k$ with the common part in $G[X]$. It follows
that $\sigma(G)\leq 3$.  $\Box$

Each $X$-component $G_i$ is smaller than $G$. In an inductive approach, we can decide
whether $\sigma(G_i)\leq 3$ for all $X$-components $G_i$ and obtain the corresponding
optimal spanning tree $T_i$ for $1\leq i\leq k$. The crucial issue is to make all $T_i$
coincident in $G[X]$. To do this, we try to find the inner star-set $S$ satisfying (II)
as the common part of all $T_i$ inside $G[X]$.

We first consider a simple case of this inductive approach.  We may encounter a neighbor
cut $X=N_G[v]$ which decomposes $G$ into $G_1$ and $G_2$ so that $\sigma(G_1)\leq 3$ and the
optimal spanning tree $T_1$ is a double star with central edge $uv$ in $G[X]$. Without loss
of generality we may assume that $T_1$ consists of a spanning star of $G[X]$ with center
$v$ and a star of $\{u\}\cup V_1$ with center $u$. Then we can ignore $G_1$ and the decision
problem of $\sigma(G)\leq 3$ can be reduced to only determining $\sigma(G_2)\leq 3$. In fact,
if $\sigma(G_2)\leq 3$ and an optimal spanning tree $T_2$ is obtained in $G_2$, then we can
construct another optimal spanning tree $T'_2$ of $G_2$ by replacing $T_2\cap G[X]$ to
$T_1\cap G[X]$ (the spanning star of $G[X]$). Thus $T_1$ and $T'_2$ are coincident in $G[X]$
and $T=T_1\cup T'_2$ is an optimal spanning tree of $G$ with $\sigma(G)\leq 3$. If $\sigma(G_2)>3$,
then $\sigma(G)>3$ by Lemma 4.3. This simplification is called a {\it reduction}. Since all
possible reductions can be found in searching for vertex cuts in $O(n^2)$ time, we assume that
the reduction computations have been done beforehand.

If a neighbor cut $X=N_G[v]$ contains an inner star-set $S=N_T[v]$ of an optimal spanning
tree $T$ (namely, $X$ comes from the inner star-set $S$ in Lemma 4.2), then it is called a
{\it normal neighbor cut}. If $\sigma(G)\leq 3$, then there exists a neighbor cut which is
normal (by Lemma 4.2). Naturally, it is possible that a neighbor cut is not normal, as
$S=N_T[v]$ may be an outer star-set of the optimal spanning tree $T$. However, we have
the following useful property.

{\bf Lemma 4.4.} \  If $D_T(G)>3$, $\sigma(G)\leq 3$, and a neighbor cut $N_G[v]$ under
consideration is not normal, then there must be a neighbor $u$ of $v$ such that $N_G[u]$ is
a normal neighbor cut.

{\bf Proof.} \  Suppose that $N_T[v]$ is an outer star-set of an optimal spanning tree $T$,
and $u$ is the unique contact vertex in $N_T(v)$. If $N_T[u]$ is an inner star-set of $T$,
then $N_G[u]$ is a normal neighbor cut, as required. Otherwise $N_T[u]$ is an outer star-set
of $T$. Thus $N_T(u)$ contains the unique contact vertex $v$. Hence $T$ is a double star with
central edge $uv$, and so $D(T)=3$, contradicting $D_T(G)>3$.

Moreover, it is possible that $v$ is a pendent vertex of an optimal spanning tree $T$ and $N_G[v]$
is a neighbor cut which is not normal. Let $u$ be the unique neighbor of $v$ in $T$. If $N_T[u]$
is an inner star-set of $T$, then $N_G[u]$ is a normal neighbor cut, as required. Otherwise
$N_T[u]$ is an outer star-set of $T$. Suppose that the vertex cut $X=N_G[v]$ decomposes $G$
into $G_1:=G[X\cup V_1]$ and $G_2:=G[X\cup V_2]$ (or more). Since $N_T[u]$ is an outer star-set
and $N_G[u]$ is not a neighbor cut, we see that $X\cup V_1\subseteq N_T[u]$ or $X\cup V_2\subseteq N_T[u]$.
Consider the former. We can construct a double star $T_1$ in $G_1$ by the spanning star of $N_G[v]$
with center $v$ and the star of $\{u\}\cup V_1$ with center $u$. It follows that $\sigma(G_1)\leq 3$.
Therefore, $X=N_G[v]$ defines a reduction. However, this situation does not occur by the assumption
just before this lemma. The proof is complete. $\Box$

To decide whether $\sigma(G)\leq 3$, for a normal neighbor cut $X=N_G[v]$, it suffices to find
the inner star-set $S\subseteq X$ satisfying (II) of Theorem 3.2. If $X=N_G[v]$ is not normal,
then we try to find a neighbor $u$ of $v$ such that $N_G[u]$ is a normal neighbor cut. If there
exists no such normal neighbor cut in the end, then we assert that $\sigma(G)>3$. Overall, our
goal is to find a normal neighbor cut that contains an inner star-set satisfying Theorem 3.2.

Suppose now that $X=N_G[v]$ is a neighbor cut of $G$ and $G$ is decomposed by $X$ into $G_1,G_2,
\ldots,G_k$, where $G_i:=G[X\cup V_i]$ ($1\leq i\leq k$). We may assume that each $G_i$ has an
optimal spanning tree $T_i$ with $\sigma(G_i)\leq 3$. For otherwise we have decided that $\sigma(G)>3$.
For each $T_i$, a subtree $F_i\subseteq T_i$ is called the {\it fixed subtree} of $T_i$ if \\
{\indent} (i) $F_i\subseteq (T_i-v)\cap G[X]$;\\
{\indent} (ii) $(T_i\setminus G[X])\cup F_i$ is connected;\\
{\indent} (iii) $F_i$ is minimal with respect to (i) and (ii).

This $F_i$ is the necessary part of $T_i$ in $G[X]$ that connects the part of $T_i$ outside $G[X]$. If
$F_i$ is a single vertex $x_i$, then $x_i$ is the only vertex connecting the edges in $T_i\setminus G[X]$.
Thus it will be taken as a contact vertex in the inner star-set $S\subseteq X$. We assume below that this
case has been considered. Otherwise $F_i$ contains some edges, which are called {\it fixed edges}. Only
the fixed edges in $G[X]$ can be adjacent to the edges in $T_i\setminus G[X]$. For a required optimal
spanning tree $T$ composed of $T_1,T_2,\ldots,T_k$, the fixed edges in $F_i$ will be contained in the
common part $T\cap G[X]$.

We first concentrate our attention on the case where each $F_i$ contains exactly one edge $e_i$. Let
$F:=\{e_1,e_2,\ldots,e_r\}$ ($r\leq k$) be the set of fixed edges in $G[X]$, where $e_i=a_ib_i$ ($1\leq i\leq r$)
with all $a_i$'s and $b_i$'s being distinct. So, $F$ is simply a matching in $G[X]$ (all edges are nonadjacent).
For an optimal spanning tree $T$ containing these fixed edges, if $S=\{v\}\cup N_T(v)$ is an inner star-set cut,
then each edge $e_i=a_ib_i$ must have one end (either $a_i$ or $b_i$) contained in $S$. We denote this end vertex
by $x_i$, which is a contact vertex of $S$ and $vx_i$ is an edge in the star $T\cap G[S]$. Meanwhile, we denote
the other end vertex of $e_i$ by $y_i$. Namely, we have $\{x_i,y_i\}=\{a_i,b_i\}$ (if $x_i=a_i$ then $y_i=b_i$).
By Theorem 3.2 (as well as the equivalent statements in Proposition 3.3), if $\sigma(G)\leq 3$, then for any
two contact vertices $x_i$ and $x_j$ of $S$, $N_T(x_i)\setminus S$ and $N_T(x_j)\setminus S$ are contained
in different components of $G-S$. Note that $y_i\in N_T(x_i)\setminus S$ and $y_j\in N_T(x_j)\setminus S$.
Hence $y_iy_j\notin E(G)$ for any $i$ and $j$ with $1\leq i<j\leq r$. Thus these $y_1,y_2,\ldots,y_r$ form
an independent (stable) set in $G[X]$. Therefore, to find the required vertex cut $S\subseteq X$ containing
$\{x_1,x_2,\ldots,x_r\}$, it suffices to determine an independent set $C=\{y_1,y_2,\ldots,y_r\}$ that
covers $F$ (each edge in $F$ has one end in $C$). This leads to the following decision problem.

{\bf The Independent Cover Set (ICS) Problem.} \  We are given a set $F=\{e_i: e_i=a_ib_i,1\leq i\leq r\}$
of fixed edges which is a matching in $G[X]$. Denote $V(F):=\{a_i:1\leq i\leq r\}\cup \{b_i:1\leq i\leq r\}$.
Is there an independent set $C\subseteq V(F)$ in $G[V(F)]$ with $|C|=r$ that covers $F$?

This is a special case of the maximum independent (stable) set problem (where $|C|=r$ means that $C$ is
maximum). Clearly, the exhaustion method of taking an end from each edge in $F$ has time bound $O(2^r)$.
However, we show that it admits a polynomial-time algorithm as follows.

{\bf Lemma 4.5.} The ICS problem can be solved by an algorithm in $O(r^2)$ time, where $r=|F|$.

{\bf Proof.} By induction on $r$. The problem is trivial for $r=1$. Assume now that $r\geq 2$ and the
assertion holds for smaller $r$.

For an edge $e_i$, if there is another edge $e_j$ such that an end (say $a_i$) of $e_i$ is adjacent to both
$a_j$ and $b_j$, then this end cannot be contained in any independent set $C$. Thus this end is deleted and
the other end (say $b_i$) should be taken into $C$. In this situation, the edge $e_i$ is said to be {\it reduced}.
The proof falls into three cases as follows.

{\bf Case 1.} \ There exists a reduced edge.

Suppose that an end (say $a_i$) of $e_i$ is adjacent to both $a_j$ and $b_j$, then this end is deleted. If both $a_i$
and $b_i$ are deleted, then stop and return the answer `no' for the decision problem (no independent set can cover $e_i$).
If exactly one end of $a_i$ and $b_i$ (say $a_i$) is deleted, then the other end (say $b_i$) is taken into $C$. By
the inductive hypothesis for $F':=F\setminus\{e_i\}$, the ICS problem of $F'$ can be solved in $O((r-1)^2)$ time.
If this problem for $F'$ has answer `no', then the same answer for $F$. If $C$ is an independent set that covers $F'$,
then $C\cup \{b_i\}$ is a solution of $F$.

If there are no reduced edges in $F$, then for any two edges $e_i,e_j\in F$, each end of $e_i$ is  independent
to either end of $e_j$, and vice versa. Let $E(\{a_i,b_i\},\{a_j,b_j\})$ denote the set of edges between $\{a_i,b_i\}$
and $\{a_j,b_j\}$ in $G[V(F)]$. Then $e_i$ and $e_j$ have the following relations (see Figure 3):\\
{\indent} $\bullet$ $e_i$ and $e_j$ are {\it parallel} if $E(\{a_i,b_i\},\{a_j,b_j\})=\{a_ia_j,b_ib_j\}$
or $\{a_ib_j,b_ia_j\}$.\\
{\indent} $\bullet$ $e_i$ and $e_j$ are {\it serial} if $E(\{a_i,b_i\},\{a_j,b_j\})=\{a_ia_j\}$ or
$\{b_ib_j\}$ or $\{a_ib_j\}$ or $\{b_ia_j\}$.\\
{\indent} $\bullet$ $e_i$ and $e_j$ are {\it unrelated} if $E(\{a_i,b_i\},\{a_j,b_j\})=\emptyset$.

\begin{center}
\setlength{\unitlength}{0.4cm}
\begin{picture}(32,8)
\multiput(0,4)(4,0){2}{\circle*{0.3}}
\multiput(0,7)(4,0){2}{\circle*{0.3}}
\put(0,4){\line(1,0){4}} \put(0,7){\line(1,0){4}}
\put(0,4){\line(0,1){3}} \put(4,4){\line(0,1){3}}

\multiput(7,4)(4,0){2}{\circle*{0.3}}
\multiput(7,7)(4,0){2}{\circle*{0.3}}
\put(7,4){\line(1,0){4}} \put(7,7){\line(1,0){4}}
\put(7,4){\line(4,3){4}} \put(7,7){\line(4,-3){4}}

\multiput(14,4)(4,0){2}{\circle*{0.3}}
\multiput(14,7)(4,0){2}{\circle*{0.3}}
\put(14,4){\line(1,0){4}} \put(14,7){\line(1,0){4}}
\put(14,4){\line(0,1){3}}

\multiput(21,4)(4,0){2}{\circle*{0.3}}
\multiput(21,7)(4,0){2}{\circle*{0.3}}
\put(21,4){\line(1,0){4}} \put(21,7){\line(1,0){4}}
\put(21,7){\line(4,-3){4}}

\multiput(28,4)(4,0){2}{\circle*{0.3}}
\multiput(28,7)(4,0){2}{\circle*{0.3}}
\put(28,4){\line(1,0){4}} \put(28,7){\line(1,0){4}}

\put(-1,7.2){\makebox(1,0.5)[l]{\small $a_i$}}
\put(-1,3.3){\makebox(1,0.5)[l]{\small $a_j$}}
\put(4.3,7.2){\makebox(1,0.5)[l]{\small $b_i$}}
\put(4.3,3.3){\makebox(1,0.5)[l]{\small $b_j$}}

\put(6,7.2){\makebox(1,0.5)[l]{\small $a_i$}}
\put(6,3.3){\makebox(1,0.5)[l]{\small $a_j$}}
\put(11.3,7.2){\makebox(1,0.5)[l]{\small $b_i$}}
\put(11.3,3.3){\makebox(1,0.5)[l]{\small $b_j$}}

\put(13,7.2){\makebox(1,0.5)[l]{\small $a_i$}}
\put(13,3.3){\makebox(1,0.5)[l]{\small $a_j$}}
\put(18.3,7.2){\makebox(1,0.5)[l]{\small $b_i$}}
\put(18.3,3.3){\makebox(1,0.5)[l]{\small $b_j$}}

\put(20,7.2){\makebox(1,0.5)[l]{\small $a_i$}}
\put(20,3.3){\makebox(1,0.5)[l]{\small $a_j$}}
\put(25.3,7.2){\makebox(1,0.5)[l]{\small $b_i$}}
\put(25.3,3.3){\makebox(1,0.5)[l]{\small $b_j$}}

\put(27,7.2){\makebox(1,0.5)[l]{\small $a_i$}}
\put(27,3.3){\makebox(1,0.5)[l]{\small $a_j$}}
\put(32.3,7.2){\makebox(1,0.5)[l]{\small $b_i$}}
\put(32.3,3.3){\makebox(1,0.5)[l]{\small $b_j$}}

\put(-0.5,2.1){\makebox(1,0.5)[l]{\small (a) parallel}}
\put(6.5,2.1){\makebox(1,0.5)[l]{\small (b) parallel}}
\put(13.9,2.1){\makebox(1,0.5)[l]{\small (c) serial}}
\put(20.9,2.1){\makebox(1,0.5)[l]{\small (d) serial}}
\put(27.5,2.1){\makebox(1,0.5)[l]{\small (e) unrelated}}
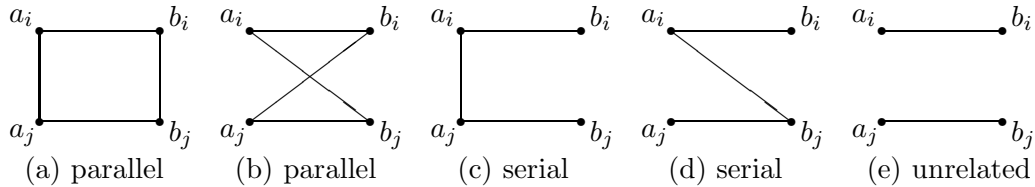
\put(4,0){\makebox(1,0.5)[l]{\small Figure 3. Relations of $e_i$ and $e_j$ for non-reduced edges.}}
\end{picture}
\end{center}

{\bf Case 2.} \ There are some parallel edges.

Suppose that $e_i$ and $e_j$ are parallel. Then $E(\{a_i,b_i\},\{a_j,b_j\})=\{a_ia_j,b_ib_j\}$ or
$\{a_ib_j,b_ia_j\}$. We first consider $E(\{a_i,b_i\},\{a_j,b_j\})=\{a_ia_j,b_ib_j\}$. Hence $\{a_i,b_j\}$
and $\{b_i,a_j\}$ are independent. If there exists an independent set $C$ that covers $F$ (thus it covers
$e_i$ and $e_j$), then $\{a_i,b_j\}\subseteq C$ or $\{b_i,a_j\}\subseteq C$. We merge $e_i$ and $e_j$ into
a single edge $uv$ by identifying $a_i$ and $b_j$ to $u$ and identifying $b_i$ and $a_j$ to $v$. Let $F'$
be the resulting edge set obtained from $F$ by this operation. Then $|F'|<|F|$. By the inductive hypothesis
for $F'$, the ICS problem of $F'$ can be solved in $O((r-1)^2)$ time. If this problem for $F'$ has answer
`no', then the same answer for $F$. Otherwise it has answer `yes' and $C$ is an independent set that covers
$F'$. If $u\in C$ (the case $v\in C$ is symmetric), then $(C\setminus \{u\})\cup \{a_i,b_j\}$ is an independent
set that covers $F$. So, the ICS problem of $F$ has answer `yes'. Similarly, if $E(\{a_i,b_i\},\{a_j,b_j\})=
\{a_ib_j,b_ia_j\}$, then $\{a_i,a_j\}$ and $\{b_i,b_j\}$ are independent. We merge $e_i$ and $e_j$ into a
single edge by identifying $a_i$ and $a_j$ and identifying $b_i$ and $b_j$. The remaining argument is the
same as before.

{\bf Case 3.} \ All edges have the serial or unrelated relations.

Suppose $e_i$ and $e_j$ are serial. Without loss of generality, assume that $E(\{a_i,b_i\},\{a_j,b_j\}) \\ =\{a_ib_j\}$
(see Figure 3 (d)). Then $\{a_i,a_j\}$ and $\{b_i,b_j\}$ are independent. We merge $e_i$ and $e_j$ into a single edge
$uv$ by identifying $a_i$ and $a_j$ to $u$ and identifying $b_i$ and $b_j$ to $v$. Let $F'$ be the resulting edge set
obtained from $F$ by this operation. By the inductive hypothesis for $F'$, the ICS problem of $F'$ can be solved in
$O((r-1)^2)$ time. If $C$ is an independent set that covers $F'$ and $u\in C$, then $(C\setminus\{u\})\cup \{a_i,a_j\}$
is the solution for $F$. Moreover, if $e_i$ is an edge unrelated to all other edges in $F$, then let $F':=F\setminus\{e_i\}$.
If $C$ is the independent set that covers $F'$, then $C$ plus any end of $e_i$ is a solution of $F$. This completes the
proof. $\Box$

We further generalize Lemma 4.5 to the case where the fixed edges are not isolated. For the fixed subtree $F_i$
of $T_i$ with $1\leq i\leq k$, let $F:=\bigcup_{1\leq i\leq k}F_i$. For an optimal spanning tree $T$ containing
$F$, if $S=\{v\}\cup N_T(v)$ is an inner star-set, then for each fixed subtree $F_i$ of $F$, there must
be a vertex $x_i\in V(F_i)$ which is a contact vertex in $S$ and $vx_i$ is an edge in the star $T\cap G[S]$.
By Theorem 3.2, if $\sigma(G)\leq 3$, then for any two contact vertices $x_i$ and $x_j$ of $S$, $N_T(x_i)
\setminus S$ and $N_T(x_j)\setminus S$ are contained in different components of $G-S$, that is, $V(F_i)
\setminus\{x_i\}$ and $V(F_j)\setminus\{x_j\}$ are contained in different $S$-branches (any vertex in
$V(F_i)\setminus\{x_i\}$ and any vertex in $V(F_j)\setminus\{x_j\}$ are independent). This gives rise
to the following decision problem.

{\bf The Consistent Fixed Edges (CFE) Problem.} Suppose that $X=N_G[v]$ is a neighbor cut of $G$ and $F_i$
is the set of fixed edges with respect to the optimal spanning tree $T_i$ of $G_i:=G[X\cup V_i]$ ($1\leq i\leq k$).
Is there an $x_i\in V(F_i)$ for each $F_i$ with $1\leq i\leq k$ such that any vertex in $V(F_i)\setminus\{x_i\}$
and any vertex in $V(F_j)\setminus\{x_j\}$ are independent for $i\neq j$? If the answer is affirmative, the sets
$F_1,F_2,\ldots,F_k$ are said to be {\it consistent}.

{\bf Lemma 4.6.} \ The Consistent Fixed Edges (CFE) Problem can be solved in $O(r^2)$ time, where $r=|F|$.

{\bf Proof.} \ We have considered the case where $|F_1|=|F_2|=\cdots=|F_k|=1$ in the ICS problem. Now we assume
that $|F_i|>1$ for some $i$ with $1\leq i\leq k$. There are three cases to consider:\\
{\indent} {\bf Case 1.} There is an $F_i$ which has radius at least three. For any $x_i\in V(F_i)$, there must a
vertex $y\in V(F_i)$ such that $d_{F_i}(x_i,y)\geq 3$. Since $vy\in E(G)$, we have $d_T(v,y)\geq 4$ for any
optimal spanning tree $T$, which contradicts that $\sigma(G)\leq 3$. Hence we answer `no' for the decision problem. \\
{\indent} {\bf Case 2.} There is an $F_i$ which has radius two. If $x_i$ is not the center of $F_i$, then there
is $y\in V(F_i)$ such that $d_{F_i}(x_i,y)\geq 3$. Thus we answer `no' as the previous case. Otherwise
$d_{F_i}(x_i,y)=2$ for the center $x_i$ and a peripheral vertex $y$ of $F_i$. If $y$ is adjacent to a vertex
$z\in V(F_j)$ for $j\neq i$, then we have $d_T(y,z)\geq 4$ for any optimal spanning tree $T$, which contradicts
that $\sigma(G)\leq 3$. Hence we answer `no' for the CFE problem. Otherwise $y$ is independent to all vertices
in $V(F_j)$ for $j\neq i$. So, we may delete $y$ from $V(F_i)$, as it does not matter whether it is incorporated
in $F_i$. Thus the decision is reduced to the next case where each $F_i$ is a star. \\
{\indent} {\bf Case 3.} All $F_1,F_2,\ldots,F_k$ are stars (including single-edge stars). Let $\{F_i: i\in I\}$
be the set of all subtrees $F_i$ with $|F_i|=1$. We first consider the ICS problem on $\{F_i: i\in I\}$. If its
answer is `no', then we return the answer `no' for the CFE problem. Otherwise (its answer is `yes') we have determined
the contact vertex $x_i$ for each $F_i$ for $i\in I$. Now we consider the problem for all stars $F_i$ with $|F_i|>1$
($i\notin I$). For each star $F_i$ with $|F_i|>1$, the contact vertex $x_i$ must be at the center of the star. For
otherwise $d_{F_i}(x_i,y)=2$ for two leaves $x_i$ and $y$, which has been discussed in the previous case. Hence
the contact vertices $x_i$ of all $F_i$ ($1\leq i\leq k$) are determined. In this situation, the CFE problem is
affirmative if and only if for any $i$ and $j$ with $1\leq i<j\leq k$, the vertices in $V(F_i)\setminus\{x_i\}$
and the vertices in $V(F_j)\setminus\{x_j\}$ are independent. This condition can be decided for all pairs of
fixed edges $x_iy$ and $x_jz$ with $i\neq j$ whether $y$ and $z$ are independent. This takes $O({r\choose 2})
=O(r^2)$ time.

To summarize, the CFE problem can be solved in $O(r^2)$ time, completing the proof. $\Box$

If $X=N_G[v]$ is a normal neighbor cut of $G$ and the CFE problem has answer `yes', then we can construct an
optimal spanning tree $T$ containing all consistent subtrees $F_1,F_2,\ldots,F_k$ in $G[X]$. Hence we assert
that $\sigma(G)\leq 3$. If the CFE problem has answer `no', that is, $F_1,F_2,\ldots,F_k$ are not consistent,
then we have two possibilities: either $\sigma(G)>3$ or $X$ is not a normal neighbor cut ($N_T[v]$ may be an
outer star-set of the optimal spanning tree $T$). In the light of Lemma 4.4, we can find a neighbor $u$ of $v$
such that $X'=N_G[u]$ is a normal neighbor cut. If this possibility occurs, we need a little more work to do.

Recall that in the decision algorithm of Theorem 2.1 for $\sigma(G)=2$, the set of 3-connected components
is decomposed by the 2-vertex cuts into a tree structure, called the {\it decomposition tree} (see
\cite{Bondy89,Bondy08,Cai95}). Now, we follow this approach. Let $X=N_G[v]$ be a neighbor cut and $G$ is
decomposed into $X$-components $G_i:=G[X\cup V_i]$ ($1\leq i\leq k$). Then $G$ is taken as the root of the
tree and $G_1,G_2,\ldots,G_k$ are the children of $G$. To decide whether $\sigma(G_i)\leq 3$ for a $G_i$,
we can find a neighbor cut in $G_i$ and carry out the decomposition to create the children of $G_i$. This
inductive approach can be performed until we reach an $X$-component which contains no neighbor cuts. In
this way, we obtain the decomposition tree, denoted by ${\mathscr T}$. For this directed tree ${\mathscr T}$,
each node is given a {\it level} below: the root has {\it level 0}; for a node of {\it level $i$}, all
its children have {\it level $i+1$}. A decomposition algorithm proceeds on this tree ${\mathscr T}$
as follows.

{\noindent}{\bf Stretch-three Recognition Algorithm} \\
\hspace*{0.2cm}\line(1,0){180}\\
{\noindent}(1) We store $G$ into ${\mathscr T}$ as a single vertex (the root).\\
{\noindent}(2) Let $H$ be the last node stored in ${\mathscr T}$. Perform the Basic Procedure for $H$. If $D_T(H)\leq 3$,\\
{\indent} then we have $\sigma(H)\leq 3$ and declare that $H$ is a {\it terminal} of ${\mathscr T}$ (stop growing). Other-\\
{\indent} wise $D_T(H)>3$. \\
{\noindent}(3) Search for a neighbor cut $X=N_H[v]$ in $H$. If no such neighbor cut, then we return \\
{\indent} answer `no' for the decision problem. Otherwise, let $X=N_H[v]$ be a neighbor cut \\
{\indent} of $H$ which decomposes $H$ into $X$-components $H_1,H_2,\ldots,H_k$. Then $H_1,H_2,\ldots,H_k$ \\
{\indent} are the children of $H$ and stored into ${\mathscr T}$.\\
{\noindent}(4) If all leaves of ${\mathscr T}$ are terminal, then go next. Otherwise go back to Step (2).\\
{\noindent}(5) If $|{\mathscr T}|=1$, then we return answer `yes' for the decision problem. \\
{\noindent}(6) Let $H$ be a node in ${\mathscr T}$ which have children $H_1,H_2,\ldots,H_k$ with maximum level. Then\\
{\indent} these children $H_i$ are terminal leaves and $\sigma (H_i)\leq 3$ with optimal spanning tree $T_i$.\\
{\indent} We construct an optimal spanning tree of $H$ by solving the CFE problem of $T_1,T_2,$ \\
{\indent} $\ldots,T_k$.\\
{\noindent}(7) If the answer of the CFE problem is `yes', then we obtain an optimal spanning tree \\
{\indent} of $H$ with $\sigma(H)\leq 3$. We change $H$ to be a terminal leaf and delete the leaves $H_1,H_2,$\\
{\indent} $\ldots,H_k$ from ${\mathscr T}$. Go back to (5). Otherwise, if the answer of the CFE problem is `no', \\
{\indent} then we return answer `no' for the decision problem. \\
\hspace*{0.2cm}\line(1,0){180}\\

Henceforth, we denote by $\eta$ the total number of $X$-components during the decomposition process. Since there are at
most $n$ neighbor cuts and for each neighbor cut $X$, there are $k$ $X$-components ($k\leq n$), it follows that $\eta=O(n^2)$.
So, the $O(\eta n^3)$-time algorithm is indeed an $O(n^5)$-time algorithm.

{\bf Theorem 4.7.} \ For any 2-connected graph $G$, deciding whether $\sigma(G)\leq 3$ can be performed by the
Stretch-three Recognition Algorithm in $O(\eta n^3)$ time, where $\eta=O(n^2)$ is the total number of $X$-components
in the algorithm.

{\bf Proof.} \  By induction on $\eta$. When $\eta=1$, $G$ contains no neighbor cuts. By Lemmas 4.1 and 4.2,
$\sigma(G)\leq 3$ can be decided by $D_T(G)\leq 3$ in $O(mn)$ time. So the assertion holds for $\eta=1$.
Assume now that $\eta>1$ and the algorithm has been established for smaller value of $\eta$.

For the first stage (Steps (1)-(4)), we construct the decomposition tree ${\mathscr T}$. For the second stage
(Steps (5)-(7)), we perform the backtrack operations from leaves to root on the tree ${\mathscr T}$. We begin
with the Basic Procedure for the root $G$. If $D_T(G)\leq 3$, then $\sigma(G)\leq 3$ and the decision is finished.
Otherwise $D_T(G)>3$. If there exists no neighbor cut in $G$, then we have $\sigma (G)>3$ by Lemma 4.2 and the
decision is also finished. Otherwise $G$ is decomposed by a neighbor cut $X=N_G[v]$ into $X$-components $G_1,G_2,
\ldots,G_k$. By Lemma 4.3, the decision process is reduced to deciding $\sigma(G_i)\leq 3$ for $1\leq i\leq k$ by the
same method. If we found one $\sigma(G_i)>3$, then the decision is also finished by answer `no'. Otherwise we obtain
all optimal panning trees $T_i$ for $\sigma(G_i)\leq 3$. If all $T_i$ are coincident in $G[X]$, then we end up with
a spanning tree $T$ of $G$ by connecting them together (by Lemma 4.3). Otherwise, we define the fixed subtree $F_i$
in $T_i$ ($1\leq i\leq k$) and solve the Consistent Fixed Edges Problem. If the subtrees $F_1,F_2,\ldots,F_k$ are
consistent, then we assert that $\sigma(G)\leq 3$. Otherwise, $F_1,F_2,\ldots,F_k$ are not consistent. Then either
$\sigma(G)>3$ or the neighbor cut $X=N_G[v]$ is not normal. For the latter possibility, it follows from Lemma 4.4
that there is a neighbor $u$ of $v$ such that $N_G[u]$ is a normal neighbor cut. However, in Step (6), the children
$H_1,H_2,\ldots,H_k$ are terminal leaves with $\sigma (H_i)\leq 3$. They contain no such neighbor cuts. So the only
possibility is $\sigma(G)>3$ and the decision is finished. In summary, we perform the Basic Procedure in $O(mn)
=O(n^3)$ time and solve the Consistent Fixed Edges Problem in $O(n^2)$ time for each node of ${\mathscr T}$.

Finally we consider the runtime of the algorithm. Suppose that the algorithm does not stop at the first step. Then we
take a neighbor cut $X=N_G[v]$ of $G$ and $G$ is decomposed by $X$ into $G_1,G_2,\ldots,G_k$. Let $n_i$ be the order
of $G_i$ and $\eta_i$ the total number of $X$-components in $G_i$ ($1\leq i\leq k$). Then $\eta_1+\eta_2+\cdots+\eta_k=\eta$.
It is clear that $n_i<n$ for all $i$ with $1\leq i\leq k$. By the inductive hypothesis, the algorithm for deciding whether
$\sigma(G_i) \leq 3$ has been established, whose complexity is $O(\eta_in_i^3)$. The overall complexity for all $G_i$ with
$1\leq i\leq k$ is $\sum_{1\leq i\leq k} O(\eta_in_i^3)$. By $\sum_{1\leq i\leq k}\eta_in_i^3\leq (\sum_{1\leq i\leq k}
\eta_i)n^3=\eta n^3$, we have the time bound $O(\eta n^3)$. This completes the proof. $\Box$

Let us see an example in Figure 4, which is due to \cite{Panda07,Panda09} as a counter example for
directed path graphs admitting $\sigma(G)\leq 3$. In this graph $G$, $N_G[v_i]$ ($1\leq i\leq 6$)
are neighbor cuts. By symmetry, we only consider $X=N_G[v_1]=\{v_1,a,b,c,v_2,v_3,v_4,v_5,v_6\}$ (the
spanning star of $X$ is shown by heavy lines in Figure 4). Then the vertex sets of components of $G-X$
are $V_1=\{d\}$, $V_2=\{e,f,g,h\}$, and $V_3=\{i,j,k,l\}$. Hence the $X$-components are as follows:
$$ G_1:=G[X\cup\{d\}], \quad G_2:=G[X\cup\{e,f,g,h\}],\quad G_3:=G[X\cup\{i,j,k,l\}].$$
Every $G_i$ has an optimal spanning tree $T_i$ with $\sigma(G_i,T_i)=3$ ($1\leq i\leq 3$) as follows:
$$ T_1=\{v_1v_2,v_1v_3,v_1v_4,v_1v_5,v_1v_6,v_1a,v_1b,v_1c,v_2d\},$$
$$ T_2=\{v_3v_2,v_3v_1,v_3v_6,v_3v_5,v_3v_4,v_3e,v_3f,v_4g,v_4h,v_1a,v_1b,v_1c\},$$
$$ T_3=\{v_6v_1,v_6v_2,v_6v_3,v_6v_4,v_6v_5,v_5i,v_5j,v_6k,v_6l,v_1a,v_1b,v_1c\}.$$

However, $T_1,T_2,T_3$ are not coincident in $G[X]$. We consider the fixed subtree $F_i$ of $T_i$
for $i=1,2,3$. Observe that $V_1=\{d\}$ and $F_1=\{v_2\}$ is a single vertex (which is a contact
vertex $x_1$). Moreover, $F_2$ contains a fixed edge $v_3v_4$, and $F_3$ contains a fixed edge
$v_5v_6$. Apart from $x_1=v_2$ incident with the center $v_1$, let $F=\{v_3v_4,v_5v_6\}$. If $G$
admits $\sigma(G)\leq 3$, then the ICS problem on $F$ must have answer `yes'. That is, there is an
independent set $C\subseteq V(F)$ with $|C|=2$ that covers $F$ (either $v_3$ or $v_4$ is independent
to either $v_5$ or $v_6$). This is impossible, as $\{v_3,v_4,v_5,v_6\}$ induces a clique in $G$.
Therefore $G$ does not admit $\sigma(G)\leq 3$.

\begin{center}
\setlength{\unitlength}{0.3cm}
\begin{picture}(20,20)
\multiput(8.5,3)(3,0){2}{\circle*{0.3}}
\multiput(10,5)(7.5,0){2}{\circle*{0.3}}
\multiput(0,8)(4,0){2}{\circle*{0.3}}
\multiput(16,8)(3,0){2}{\circle*{0.3}}
\multiput(0,14)(4,0){2}{\circle*{0.3}}
\multiput(16,14)(3,0){2}{\circle*{0.3}}
\multiput(2,17)(8,0){2}{\circle*{0.3}}
\multiput(8.5,19)(3,0){2}{\circle*{0.3}}
\multiput(17.5,16)(0,0){1}{\circle*{0.3}}
\multiput(2,5)(0,0){1}{\circle*{0.3}}

\put(0,8){\line(1,0){4}} \put(0,14){\line(1,0){4}}
\put(0,8){\line(2,-3){2}} \put(2,5){\line(2,3){2}}
\put(0,14){\line(2,3){2}} \put(4,14){\line(-2,3){2}}
\put(0,8){\line(2,3){4}} \put(4,8){\line(-2,3){4}}
\put(0,8){\line(0,1){6}} \put(4,8){\line(0,1){6}}

\put(10,5){\line(-2,1){6}} \put(10,5){\line(2,1){6}}
\put(10,17){\line(-2,-1){6}} \put(10,17){\line(2,-1){6}}
\put(10,5){\line(-2,3){6}} \put(10,5){\line(2,3){6}}
\put(10,17){\line(-2,-3){6}} \put(10,17){\line(2,-3){6}}
\put(4,8){\line(1,0){12}} \put(4,14){\line(1,0){12}}
\put(16,8){\line(0,1){6}} \put(10,5){\line(0,1){12}}
\put(4,8){\line(2,1){12}} \put(4,14){\line(2,-1){12}}

\put(10,17){\line(-3,4){1.5}} \put(10,17){\line(3,4){1.5}}
\put(8.5,19){\line(1,0){3}} \put(16,14){\line(1,0){3}}
\put(16,14){\line(3,4){1.5}} \put(19,14){\line(-3,4){1.5}}
\put(11.5,19){\line(2,-1){6}}
\qbezier(10,17)(13.75,16.5)(17.5,16)
\qbezier(11.5,19)(13.75,16.5)(16,14)

\put(10,5){\line(-3,-4){1.5}} \put(10,5){\line(3,-4){1.5}}
\put(8.5,3){\line(1,0){3}} \put(16,8){\line(1,0){3}}
\put(16,8){\line(1,-2){1.5}} \put(10,5){\line(1,0){7.5}}
\put(17.5,5){\line(1,2){1.5}}
\qbezier(11.5,3)(13.75,5.5)(16,8)
\qbezier(11.5,3)(14.5,4)(17.5,5)

\put(8.3,17){\makebox(1,0.5)[l]{\small $v_1$}}
\put(16.4,13.1){\makebox(1,0.5)[l]{\small $v_2$}}
\put(16.4,8.3){\makebox(1,0.5)[l]{\small $v_3$}}
\put(8.3,4.7){\makebox(1,0.5)[l]{\small $v_4$}}
\put(3.7,6.8){\makebox(1,0.5)[l]{\small $v_5$}}
\put(3.7,14.8){\makebox(1,0.5)[l]{\small $v_6$}}
\put(8,19.5){\makebox(1,0.5)[l]{\small $a$}}
\put(11.4,19.5){\makebox(1,0.5)[l]{\small $b$}}
\put(17.9,16.3){\makebox(1,0.5)[l]{\small $c$}}
\put(19.4,13.8){\makebox(1,0.5)[l]{\small $d$}}
\put(19.4,7.8){\makebox(1,0.5)[l]{\small $e$}}
\put(17.8,4.6){\makebox(1,0.5)[l]{\small $f$}}
\put(11.2,2){\makebox(1,0.5)[l]{\small $g$}}
\put(8.2,2){\makebox(1,0.5)[l]{\small $h$}}
\put(1.4,4.2){\makebox(1,0.5)[l]{\small $i$}}
\put(-1,8){\makebox(1,0.5)[l]{\small $j$}}
\put(-1,14){\makebox(1,0.5)[l]{\small $k$}}
\put(1.5,17.5){\makebox(1,0.5)[l]{\small $l$}}

\thicklines{
\put(10,17){\line(-2,-1){6}} \put(10,17){\line(2,-1){6}}
\put(10,17){\line(-2,-3){6}} \put(10,17){\line(2,-3){6}}
\put(10,5){\line(0,1){12}} \put(10,17){\line(-3,4){1.5}}
\put(10,17){\line(3,4){1.5}}
\qbezier(10,17)(13.75,16.5)(17.5,16)}

\put(-3.8,0){\makebox(1,0.5)[l]{\small Figure 4. A graph $G$ not admitting $\sigma(G)\leq 3$.}}
\end{picture}
\end{center}

\begin{center}
\setlength{\unitlength}{0.3cm}
\begin{picture}(20,20)
\multiput(8.5,3)(3,0){2}{\circle*{0.3}}
\multiput(10,5)(7.5,0){2}{\circle*{0.3}}
\multiput(0,8)(4,0){2}{\circle*{0.3}}
\multiput(16,8)(3,0){2}{\circle*{0.3}}
\multiput(0,14)(4,0){2}{\circle*{0.3}}
\multiput(16,14)(3,0){2}{\circle*{0.3}}
\multiput(2,17)(8,0){2}{\circle*{0.3}}
\multiput(8.5,19)(3,0){2}{\circle*{0.3}}
\multiput(17.5,16)(0,0){1}{\circle*{0.3}}
\multiput(2,5)(0,0){1}{\circle*{0.3}}

\put(0,8){\line(1,0){4}} \put(0,14){\line(1,0){4}}
\put(0,8){\line(2,-3){2}} \put(2,5){\line(2,3){2}}
\put(0,14){\line(2,3){2}} \put(4,14){\line(-2,3){2}}
\put(0,8){\line(2,3){4}} \put(4,8){\line(-2,3){4}}
\put(0,8){\line(0,1){6}} \put(4,8){\line(0,1){6}}

\put(10,5){\line(2,1){6}}
\put(10,17){\line(-2,-1){6}} \put(10,17){\line(2,-1){6}}
\put(10,5){\line(-2,3){6}} \put(10,5){\line(2,3){6}}
\put(10,17){\line(-2,-3){6}} \put(10,17){\line(2,-3){6}}
\put(4,8){\line(1,0){12}} \put(4,14){\line(1,0){12}}
\put(16,8){\line(0,1){6}} \put(10,5){\line(0,1){12}}
\put(4,8){\line(2,1){12}} \put(4,14){\line(2,-1){12}}

\put(10,17){\line(-3,4){1.5}} \put(10,17){\line(3,4){1.5}}
\put(8.5,19){\line(1,0){3}} \put(16,14){\line(1,0){3}}
\put(16,14){\line(3,4){1.5}} \put(19,14){\line(-3,4){1.5}}
\put(11.5,19){\line(2,-1){6}}
\qbezier(10,17)(13.75,16.5)(17.5,16)
\qbezier(11.5,19)(13.75,16.5)(16,14)

\put(10,5){\line(-3,-4){1.5}} \put(10,5){\line(3,-4){1.5}}
\put(8.5,3){\line(1,0){3}} \put(16,8){\line(1,0){3}}
\put(16,8){\line(1,-2){1.5}} \put(10,5){\line(1,0){7.5}}
\put(17.5,5){\line(1,2){1.5}}
\qbezier(11.5,3)(13.75,5.5)(16,8)
\qbezier(11.5,3)(14.5,4)(17.5,5)

\put(8.3,17){\makebox(1,0.5)[l]{\small $v_1$}}
\put(16.4,13.1){\makebox(1,0.5)[l]{\small $v_2$}}
\put(16.4,8.3){\makebox(1,0.5)[l]{\small $v_3$}}
\put(8.3,4.7){\makebox(1,0.5)[l]{\small $v_4$}}
\put(3.7,6.8){\makebox(1,0.5)[l]{\small $v_5$}}
\put(3.7,14.8){\makebox(1,0.5)[l]{\small $v_6$}}
\put(8,19.5){\makebox(1,0.5)[l]{\small $a$}}
\put(11.4,19.5){\makebox(1,0.5)[l]{\small $b$}}
\put(17.9,16.3){\makebox(1,0.5)[l]{\small $c$}}
\put(19.4,13.8){\makebox(1,0.5)[l]{\small $d$}}
\put(19.4,7.8){\makebox(1,0.5)[l]{\small $e$}}
\put(17.8,4.6){\makebox(1,0.5)[l]{\small $f$}}
\put(11.2,2){\makebox(1,0.5)[l]{\small $g$}}
\put(8.2,2){\makebox(1,0.5)[l]{\small $h$}}
\put(1.4,4.2){\makebox(1,0.5)[l]{\small $i$}}
\put(-1,8){\makebox(1,0.5)[l]{\small $j$}}
\put(-1,14){\makebox(1,0.5)[l]{\small $k$}}
\put(1.5,17.5){\makebox(1,0.5)[l]{\small $l$}}

\thicklines{
\put(10,17){\line(-2,-1){6}} \put(10,17){\line(2,-1){6}}
\put(10,17){\line(2,-3){6}} \put(10,17){\line(-3,4){1.5}}
\put(10,17){\line(3,4){1.5}} \put(10,17){\line(2,-1){6}}
\put(16,14){\line(1,0){3}} \put(4,8){\line(0,1){6}}
\put(2,5){\line(2,3){2}} \put(0,8){\line(1,0){4}}
\put(0,14){\line(1,0){4}} \put(4,14){\line(-2,3){2}}
\put(10,5){\line(-3,-4){1.5}} \put(10,5){\line(3,-4){1.5}}
\put(10,5){\line(2,1){6}}\put(16,8){\line(1,0){3}}
\put(16,8){\line(1,-2){1.5}}}

\put(-3,0){\makebox(1,0.5)[l]{\small Figure 5. A graph $G$ admitting $\sigma(G)\leq 3$.}}
\end{picture}
\end{center}

By the same argument, for the graph $G':=G-v_4v_5$ shown in Figure 5, the ICS problem on
$F=\{v_3v_4,v_5v_6\}$ has answer `yes', as $C=\{v_4,v_5\}$ is independent. Hence $\sigma(G')\leq 3$
and the optimal spanning tree $T$ is shown by the heavy lines in Figure 5.

\section{Concluding remarks}
\hspace*{0.5cm} Up to now, we see that determining $\sigma(G)\leq 3$ is
polynomially solvable, while determining $\sigma(G)\leq k$ for $k\geq 4$
is NP-complete. We fill the gap of complexity.

For the extremal graph problem with respect to $\sigma(G)\geq g(G)-1$ (where
$g(G)$ is the girth of $G$), Bondy \cite{Bondy89} and Cai and Corneil \cite{Cai95}
presented the characterization for $g(G)=3$ in Theorem 2.1. For $g(G)=4$ (for
example, $G$ is bipartite), it is interesting to characterize the graphs which
have a spanning tree $T$ with every fundamental cycle being a quadrilateral. The
answer is included in Theorem 3.2, but it should be stated more precisely. For
$g(G)=5$, the Petersen graph is one of the extremal graphs. How about the others?
The graphs admitting $\sigma(G)\leq 3$ may be called the {\it star-separating
graphs}. This is a big family, including the trigraphs in \cite{Bondy89}, interval
graphs, split graphs, permutation graphs, convex graphs, etc. It is worthwhile to
study more structural features of this family.

\section*{Acknowledgments}
This work was supported by National Key R\&D Program of China under grant 2019YFB2101604.

\end{document}